\documentclass[letterpaper, 10 pt, conference]{ieeeconf}

\usepackage{epsfig}
\usepackage{amsmath}
\usepackage{amssymb} 

\title{\LARGE {\bf Natural frames and interacting particles in three dimensions
}}

\author{ \parbox{3 in}{\centering E. W. Justh \\
         Institute for Systems Research \\
         University of Maryland \\
         College Park, MD 20742, USA \\
         {\tt\small justh@umd.edu}}
         \hspace*{ 0.5 in}
         \parbox{3 in}{ \centering P. S. Krishnaprasad \\
         Institute for Systems Research and \\
         Dept. of Electrical and Computer Engineering \\
         University of Maryland \\
         College Park, MD 20742, USA \\
         {\tt\small krishna@umd.edu}}
}

\begin{document}

\maketitle
\thispagestyle{empty}
\pagestyle{empty}

\begin{abstract}

Motivated by the problem of formation control for vehicles moving 
at unit speed in three-dimensional space, we are led to models of
gyroscopically interacting particles, which require the 
machinery of curves and frames to describe and analyze. 
A Lie group formulation arises naturally, and we discuss the
general problem of determining (relative) equilibria
for arbitrary \boldmath$ G $\unboldmath-invariant controls
(where \boldmath$ G=SE(3) $ \unboldmath is a symmetry
group for the control law).
We then present global convergence (and non-collision) results
for specific two-vehicle interaction laws in three dimensions,
which lead to specific formations (i.e., relative equilibria).
Generalizations of the interaction laws to \boldmath$ n $ \unboldmath
vehicles is also discussed, and simulation results presented.

\end{abstract}

\section{Introduction}

This work is motivated by the problem of multi-vehicle formation
(or swarm) control, e.g., for meter-scale UAVs (unmanned aerial
vehicles), and builds on our earlier work on planar formation
control laws \cite{scltechrep,scl02,cdc03} by extending the key
results to the three-dimensional setting.  Some objectives of our
formation control laws are to avoid collisions between vehicles,
maintain cohesiveness of the formation, be robust to loss of 
individuals, and scale favorably to large swarms.  

In considering the problem of multi-vehicle formation control,
there is special
significance, both practically and theoretically, to modeling the
vehicles as point particles moving at a common (constant) speed.
In the language of mechanics, the individual particles are subject
to {\it gyroscopic} forces; i.e., forces which alter the direction of
motion of the particles, but which leave their speed (and hence their
kinetic energy) unchanged.  A formation control law is then a
feedback law which specifies these gyroscopic forces based on the
positions and directions of motion of the particles.
In the planar setting, gyroscopic forces serve as steering controls
\cite{scl02}.  
For particles moving in three dimensional space, we need to introduce
the notion of {\it framing of curves} to describe the effects of
gyroscopic forces on particle motion \cite{bishop,calini}.

Recently, a growing literature has emerged on planar
formation control for unit-speed vehicles, using tools from dynamical
systems theory (including pursuit models \cite{francis} and phase-coupled 
oscillator models \cite{sepulchre}), as well as graph-theoretic methods 
\cite{jadbabaie}.  An early (discrete-time) unit-speed model for
biological flocking behavior is the Vicsek model \cite{vicsek}.
Interacting particle models similar to those described in this paper 
have also found application in obstacle 
avoidance and boundary following \cite{zhang}.

\section{Curves and moving frames}

A single particle moving in three dimensional space  
traces out a trajectory 
$ \mbox{\boldmath$\gamma$\unboldmath}: [0,\infty)
 \rightarrow \mathbb{R}^3 $, 
which we assume to be at least twice continuously differentiable, 
satisfying 
$ |\mbox{\boldmath$\gamma$\unboldmath}'(s)| = 1, $ $ \forall s; $
i.e., $ s $ is the arc-length parameter of the curve (and the prime
denotes differentiation with respect to $ s $).
The direction of motion of the particle at $ s $ is the 
unit tangent vector to the trajectory, 
$ {\bf T}(s) = \mbox{\boldmath$\gamma$\unboldmath}'(s) $.
If we further restrict the speed of particle motion to be unit speed, 
then the arclength parameter $ s $ is equivalent to time $ t $,
and $ {\bf T}(t) = \dot{\mbox{\boldmath$\gamma$\unboldmath}}(t) $.  
The gyroscopic 
force vector always lies in the plane perpendicular to $ {\bf T} $, so to 
describe the effects of this force, we are compelled to introduce orthonormal 
unit vectors which span this {\it normal plane}.  Taken together with 
$ {\bf T} $, these unit vectors constitute a {\it framing} of the curve 
\boldmath$\gamma$ \unboldmath representing the particle trajectory.

\begin{figure}
\epsfxsize=9cm
\epsfbox{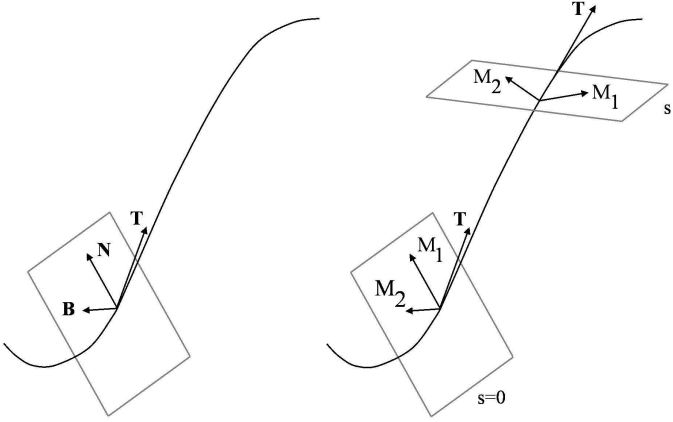}
\caption{\label{frame_3d_fig} The Frenet-Serret frame (left), and
natural Frenet frame (right), illustrated for a three-dimensional curve.}
\end{figure}

There are different framings one can choose, as is best illustrated by 
examples (see figure \ref{frame_3d_fig}).
For a curve \boldmath$\gamma$\unboldmath $(s) $ which is three times 
continuously differentiable,
and for which \boldmath$\gamma$\unboldmath ${}''(s) \ne 0 $ for all 
$ s $, the Frenet-Serret frame 
$({\bf T},{\bf N},{\bf B})$ is uniquely defined, and satisfies
\begin{eqnarray}
\label{frenetserret}
\mbox{\boldmath$\gamma$\unboldmath}'(s) 
 \hspace{-.2cm} & = & \hspace{-.2cm} {\bf T}(s), \nonumber \\
{\bf T}'(s) \hspace{-.2cm} & = & \hspace{-.2cm} \kappa(s) {\bf N}(s), 
\nonumber \\
{\bf N}'(s) \hspace{-.2cm} & = & \hspace{-.2cm} 
 - \kappa(s) {\bf T}(s) + \tau(s) {\bf B}(s),
 \nonumber \\
{\bf B}'(s)  \hspace{-.2cm} & = & \hspace{-.2cm} -\tau(s) {\bf N}(s).
\end{eqnarray}
In (\ref{frenetserret}), 
$ {\bf N}(s) $
is the unit normal vector to the curve \boldmath$\gamma$ \unboldmath
at $ s $, and $ {\bf B}(s) $ is the unit 
binormal vector
(which completes the right-handed orthonormal frame).  
The curvature function $ \kappa $ and the torsion function $ \tau $ 
are given by expressions involving the derivatives of 
\boldmath$\gamma$\unboldmath, and 
\boldmath$\gamma$\unboldmath ${}''(s) \ne 0 $ is required for
$ \tau(s) $ to be well-defined. 

Although the Frenet-Serret frame for a curve (when it exists) has a 
special status (because it is uniquely defined by the derivatives of
the curve), it is not the only choice of frame, nor is it 
necessarily the best choice.  In particular,
the requirement that $ \mbox{\boldmath$\gamma$\unboldmath}''(s) \ne 0 $
presents serious difficulties for the interaction laws we consider
in this paper.

We therefore use an alternative framing of the curve
\boldmath$\gamma$\unboldmath, the natural Frenet frame, which is also 
referred to as the Fermi-Walker frame or Relatively Parallel Adapted Frame 
(RPAF):
\begin{eqnarray}
\label{naturalfrenet}
\mbox{\boldmath$\gamma$\unboldmath}'(s)
 \hspace{-.2cm} & = & \hspace{-.2cm} {\bf T}(s), \nonumber \\
{\bf T}'(s) \hspace{-.2cm} & = & \hspace{-.2cm} 
 k_1(s){\bf M}_1 + k_2(s) {\bf M}_2, \nonumber \\
{\bf M}_1'(s) \hspace{-.2cm} & = & \hspace{-.2cm}  - k_1(s){\bf T}(s),
 \nonumber \\
{\bf M}_2'(s) \hspace{-.2cm} & = & \hspace{-.2cm}  - k_2(s){\bf T}(s).
\end{eqnarray}
In (\ref{naturalfrenet}),  
$ {\bf M}_1(s) $ and $ {\bf M}_2(s) $ are unit normal vectors which
(along with $ {\bf T}(s) $) complete a right-handed orthonormal frame.
However, there is freedom in the choice of initial conditions 
$ {\bf M}_1(0) $ and $ {\bf M}_2(0) $; once these are
specified, the corresponding natural Frenet frame for a 
twice-continuously-differentiable curve
\boldmath$\gamma$ \unboldmath
is unique.

Both (\ref{frenetserret}) and (\ref{naturalfrenet}) can be packaged
as control systems on the Lie group $ SE(3) $, the group of rigid motions in 
three-dimensional space.
(A modern reference for control systems
on Lie groups is Jurdjevic \cite{jurdjevic}.)
Here we think of $ (\kappa, \tau) $ or the {\it natural curvatures}
$ (k_1,k_2) $ as controls,
which drive the evolution of the frame and the particle position
\boldmath$\gamma$\unboldmath.

\section{Formation model}

\begin{figure}
\hspace{.5cm}
\epsfxsize=7cm
\epsfbox{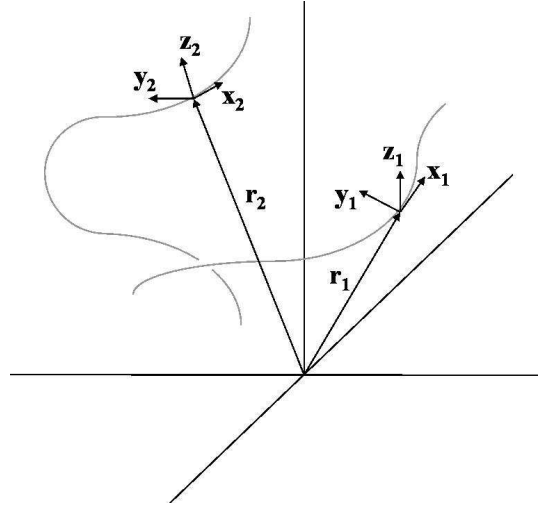}
\caption{\label{motion3dfig} Three-dimensional trajectories for two
vehicles, and their respective frames.}
\end{figure}

Figure \ref{motion3dfig} illustrates the trajectories of two
vehicles moving at unit speed, and their respective natural Frenet frames.  
The particle (i.e., vehicle) positions are denoted by $ {\bf r}_1 $
and $ {\bf r}_2 $, and the frames by $ ({\bf x}_1,{\bf y}_1,{\bf z}_1) $
and $ ({\bf x}_2,{\bf y}_2,{\bf z}_2) $, so that
\begin{eqnarray}
\label{twouavsystem3d}
\dot{\bf r}_1 = {\bf x}_1, \hspace{1.8cm}
  &&  \dot{\bf r}_2  =  {\bf x}_2,\nonumber \\
\dot{\bf x}_1 = {\bf y}_1 u_1 + {\bf z}_1 v_1, \hspace{.35cm}
  && \dot{\bf x}_2  =  {\bf y}_2 u_2 + {\bf z}_2 v_2, \nonumber \\
\dot{\bf y}_1 = -{\bf x}_1 u_1, \hspace{1.15cm}
  && \dot{\bf y}_2  = -{\bf x}_2 u_2, \nonumber \\
\dot{\bf z}_1 = -{\bf x}_1 v_1, \hspace{1.2cm}
  && \dot{\bf z}_2  = -{\bf x}_2 v_2.
\end{eqnarray}
where the controls $ (u_1,v_1) $ and $ (u_2,v_2) $ 
may be feedback functions of the position and frame variables. 

We consider control laws which depend only on relative vehicle
positions and orientations; i.e., which depend only on the {\it shape}
of the formation.  Furthermore, the {\it effect} of the controls
on each trajectory is assumed to depend only on  
$ {\bf r}_1 $, $ {\bf r}_2 $, $ {\bf x}_1 $, and $ {\bf x}_2 $,
and not on the orientation of the normal vectors within their 
respective normal planes.
 
The controls for the first vehicle can then be functions
of the relative vehicle position, $ {\bf r} = {\bf r}_2 - {\bf r}_1 $, 
the heading direction
of the second vehicle, $ {\bf x}_2 $, and the frame variables for
the first vehicle, $ ({\bf x}_1,{\bf y}_1, {\bf z}_1) $.  Thus,
\begin{eqnarray}
\label{vehctrlrestrict1}
u_1 \hspace{-.2cm} & = & \hspace{-.2cm}
 u_1({\bf r},{\bf x}_1,{\bf y}_1,{\bf z}_1,{\bf x}_2),
 \nonumber \\ 
v_1 \hspace{-.2cm} & = & \hspace{-.2cm}
 v_1({\bf r},{\bf x}_1,{\bf y}_1,{\bf z}_1,{\bf x}_2),
\end{eqnarray}
and similarly,
\begin{eqnarray}
u_2 \hspace{-.2cm} & = & \hspace{-.2cm}
 u_2({\bf r},{\bf x}_2,{\bf y}_2,{\bf z}_2,{\bf x}_1),
 \nonumber \\ 
v_2 \hspace{-.2cm} & = & \hspace{-.2cm}
 v_2({\bf r},{\bf x}_2,{\bf y}_2,{\bf z}_2,{\bf x}_1).
\end{eqnarray}
Furthermore, because the overall motion of the first vehicle 
should be independent of $ {\bf y}_1 $ and $ {\bf z}_1 $, we require
\begin{equation}
v_1({\bf r},{\bf x}_1,{\bf y}_1,{\bf z}_1,{\bf x}_2) =
 u_1({\bf r},{\bf x}_1,{\bf z}_1,-{\bf y}_1,{\bf x}_2), 
\end{equation}
and similarly,
\begin{equation}
\label{vehctrlrestrict2}
v_2({\bf r},{\bf x}_2,{\bf y}_2,{\bf z}_2,{\bf x}_1) = 
 u_2({\bf r},{\bf x}_2,{\bf z}_2,-{\bf y}_2,{\bf x}_1).
\end{equation}

Finally, we require that our control laws have a discrete
(relabling) symmetry, which corresponds to the intuitive notion
that both vehicles ``run the same algorithm.'' This implies
\begin{eqnarray}
\label{vehctrlrestrict3}
u_1(-{\bf r},{\bf x}_1,{\bf y}_1,{\bf z}_1,{\bf x}_2) 
 = u_2({\bf r},{\bf x}_2,{\bf y}_2,{\bf z}_2,{\bf x}_1), \nonumber \\
v_1(-{\bf r},{\bf x}_1,{\bf y}_1,{\bf z}_1,{\bf x}_2)
 = v_2({\bf r},{\bf x}_2,{\bf y}_2,{\bf z}_2,{\bf x}_1).
\end{eqnarray}
In this paper, the specific control laws we consider 
have the form
\begin{eqnarray}
\label{twovehiclelaw3d}
u_1 \hspace{-.2cm} & = & \hspace{-.2cm}
  F(-{\bf r},{\bf x}_1,{\bf y}_1,{\bf x}_2)
 - f(|{\bf r}|)\left(-\frac{\bf r}{|{\bf r}|}\cdot{\bf y}_1\right),
\nonumber \\
u_2 \hspace{-.2cm} & = & \hspace{-.2cm}
  F({\bf r},{\bf x}_2,{\bf y}_2,{\bf x}_1)
 - f(|{\bf r}|)\left(\frac{\bf r}{|{\bf r}|}\cdot{\bf y}_2\right),
\nonumber \\
v_1 \hspace{-.2cm} & = & \hspace{-.2cm}
  F(-{\bf r},{\bf x}_1,{\bf z}_1,{\bf x}_2)
 - f(|{\bf r}|)\left(-\frac{\bf r}{|{\bf r}|}\cdot{\bf z}_1\right),
\nonumber \\
v_2 \hspace{-.2cm} & = & \hspace{-.2cm}
  F({\bf r},{\bf x}_2,{\bf z}_2,{\bf x}_1) 
 - f(|{\bf r}|)\left(\frac{\bf r}{|{\bf r}|}\cdot{\bf z}_2\right),
\end{eqnarray}
which is a further restricted class of laws consistent with
 (\ref{vehctrlrestrict1}) - (\ref{vehctrlrestrict3}).  (We discuss
later how $ F $ and $ f $ are chosen.)

\section{Shape variables and equilibria}

The geometry of the problem of interacting particles moving at unit 
speed in the plane has been considered in earlier work 
\cite{scltechrep,scl02,cdc03}.
The unit speed constraint leads to the study of gyroscopic 
interaction forces, and the identification of the constant
kinetic energy hyper-surface with the group $ SE(2) $ of rigid motions
in the plane.  Formations or steady patterns of motion in the plane
thus become relative equilibria for particle dynamics on $ SE(2) $
\cite{scltechrep,scl02,cdc03}. 

A key difficulty in extending the above geometric perspective to
three dimensions arises from the fact that the corresponding
constant kinetic energy hyper-surface cannot be identified with
$ SE(3) $, the rigid motion group in three dimensions.  It is a 
homogeneous space $ SE(3)/SO(2) $.  However, there is considerable
advantage, particularly in the multi-particle context, to 
formulating the dynamics in terms of interacting particles in 
$ SE(3) $.

The dynamics (\ref{twouavsystem3d})
can be expressed in terms of the group variables 
$ g_1, g_2 \in G = SE(3) $ as a pair of left-invariant systems
\begin{equation}
\label{se3system}
\dot{g}_1 = g_1 \xi_1, \;\; \dot{g}_2 = g_2 \xi_2,
\end{equation}
where  $ \xi_1, \xi_2 \in {\mathfrak g}= $ the Lie algebra of $ G $.
The dynamics for $ g = g_1^{-1} g_2 $ are given by
\begin{eqnarray}
\label{gdot}
\dot{g} && \hspace{-.6cm} 
 = -g_1^{-1}\dot{g}_1 g_1^{-1} g_2 + g_1^{-1} \dot{g}_2 \nonumber \\
 && \hspace{-.6cm} = -g_1^{-1}g_1 \xi_1 g + g_1^{-1} g_2 \xi_2 \nonumber \\
 && \hspace{-.6cm} = - \xi_1 g + g \xi_2 \nonumber \\
 && \hspace{-.6cm} = g \xi,
\end{eqnarray}
where $ \xi = \xi_2 - \mbox{Ad}_{g^{-1}} \xi_1 \in {\mathfrak g} $.

Equation (\ref{gdot}), where $ \xi $ incorporates the control inputs
$ (u_1,v_1) $ and $ (u_2,v_2) $, describes the evolution of the 
{\it relative} position and {\it relative} natural Frenet frame orientation
of the pair of vehicles.  It is thus natural to consider what 
equilibria of (\ref{gdot}) exist, and then to design control laws
which stabilize those equilibria.  Equilibria of the shape dynamics
(\ref{gdot}) correspond to {\it relative equilibria} of the system
(\ref{se3system}) on $ G \times G $.  

\subsection{Shape equilibria for a two-particle system on SE(3)}

At an equilibrium shape $ g_e $ of the shape dynamics 
(\ref{gdot}), we have
\begin{equation}
\label{gexi2xi1}
g_e \xi_2(g_e) = \xi_1(g_e) g_e.
\end{equation}
To facilitate calculation, we define
\begin{eqnarray}
g_e \hspace{-.2cm} & = & \hspace{-.2cm} 
\left[ \begin{array} {c c} Q & {\bf b} \\ {\bf 0} & 1 \end{array} \right],
 \mbox{ where $ Q \in SO(3) $ and $ {\bf b} \in \mathbb{R}^3 $}, \nonumber \\
\xi_1(g_e) \hspace{-.2cm} & = & \hspace{-.2cm}
\left[ \begin{array} {c c} \hat{\Omega}_1 & {\bf e}_1 \\ {\bf 0} & 0
 \end{array} \right], \;\;
\xi_2(g_e) =
\left[ \begin{array} {c c} \hat{\Omega}_2 & {\bf e}_1 \\ {\bf 0} & 0
 \end{array} \right].
\end{eqnarray}
Then (\ref{gexi2xi1}) becomes
\begin{equation}
\label{equilblockmatrix}
\left[ \begin{array} {c c} Q & {\bf b} \\ {\bf 0} & 1 \end{array} \right]
\left[ \begin{array} {c c} \hat{\Omega}_2 & {\bf e}_1 \\ {\bf 0} & 0
 \end{array} \right] = 
\left[ \begin{array} {c c} \hat{\Omega}_1 & {\bf e}_1 \\ {\bf 0} & 0
 \end{array} \right]
\left[ \begin{array} {c c} Q & {\bf b} \\ {\bf 0} & 1 \end{array} \right],
\end{equation}
where 
$ {\bf e}_1 = \left[ \begin{array} {c c c} 1 & 0 & 0 \end{array} \right]^T $,
\begin{equation}
\label{equilblockmatrixdefn}
\Omega_1 =
  \left[ \begin{array} {c} w_1 \\ -v_1 \\ u_1 \end{array} \right], \;\;
\Omega_2 =
  \left[ \begin{array} {c} w_2 \\ -v_2 \\ u_2 \end{array} \right], 
\end{equation}
and for any 3-vector $ \Gamma = (\Gamma_1,\Gamma_2,\Gamma_3) $, 
$\hat{\Gamma} $ is the skew-symmetric matrix defined by 
\begin{equation}
\hat{\Gamma} = \left[ \begin{array} {r r r} 0 \;\; & -\Gamma_3 & \Gamma_2 \\
 \Gamma_3 & 0 \;\; & -\Gamma_1 \\ -\Gamma_2 & \Gamma_1 & 0 \;\;
 \end{array} \right].
\end{equation}
 
Note that here we allow $ \Omega_1 $ and $ \Omega_2 $ to each have the 
full three degrees of
freedom - not just the two corresponding to the natural curvatures.
The reason
for proceeding in this manner is that ultimately we recover not only the
relative equilibria of (\ref{se3system}) and (\ref{twouavsystem3d}), but 
also an interesting class
of relative {\it periodic} solutions for (\ref{twouavsystem3d}). 

From (\ref{equilblockmatrix}) we see that 
$ Q \hat{\Omega}_2 = \hat{\Omega}_1 Q $, from which it follows that
\begin{equation}
\label{rotmatrixident}
\Omega_1 = Q \Omega_2.
\end{equation}
From (\ref{equilblockmatrix}) we also obtain
$ Q {\bf e}_1 = \hat{\Omega}_1 {\bf b} + {\bf e}_1 $.  It can then be
shown that $ w_1 = w_2 $, and $ u_1^2 + v_1^2 = u_2^2 + v_2^2 $.

Introducing new variables $ w $, $ a $, $ \psi_1 $, and 
$ \psi_2 $, we can express $ \Omega_1 $ and $ \Omega_2 $ as 
\begin{equation}
\label{ctrlvectors}
\Omega_1 = \left[ \begin{array} {c} w \\ a\sin\psi_1 \\ a\cos\psi_1
 \end{array} \right], \;\;\;\;
\Omega_2 = \left[ \begin{array} {c} w \\ a\sin\psi_2 \\ a\cos\psi_2
 \end{array} \right].
\end{equation}
If (for $ a^2+w^2 \ne 0 $) we further define
\begin{equation}
\label{varphidefn}
\cos\varphi = \frac{a}{\sqrt{a^2+w^2}}, \;\;
\sin\varphi =  \frac{w}{\sqrt{a^2+w^2}},
\end{equation}
along with
\begin{eqnarray}
\label{rotmatrixdefn}
R_{\psi_j} \hspace{-.3cm} & = & \hspace{-.3cm}
 \left[ \hspace{-.15cm} \begin{array} {c c c} 1 & 0 & 0 \\ 
 0 & \cos\psi_j & -\sin\psi_j \\
 0 & \sin\psi_j & \cos\psi_j \end{array} \hspace{-.15cm} \right]
 \hspace{-.1cm}, \;
R_{\varphi} \hspace{-.1cm} = \hspace{-.1cm} 
\left[ \hspace{-.15cm}
 \begin{array} {c c c} \cos\varphi & 0 & -\sin\varphi \\ 0 & 1 & 0 \\
  \sin\varphi & 0 & \cos\varphi \end{array} \hspace{-.15cm} \right]
 \hspace{-.1cm}, 
\nonumber \\
R_{\vartheta} \hspace{-.3cm} & = & \hspace{-.3cm}
 \left[ \hspace{-.15cm} \begin{array} {c c c} 
 \cos\vartheta & -\sin\vartheta & 0 \\ 
 \sin\vartheta & \cos\vartheta & 0 \\
 0 & 0 & 1 \end{array} \hspace{-.15cm} \right],
\end{eqnarray}
where $ \vartheta \in [0,2\pi) $ is arbitrary, we see that 
(\ref{ctrlvectors}) becomes
\begin{equation}
\Omega_j = \sqrt{a^2+w^2}\; R_{\psi_j}^T R_{\varphi}^T {\bf e}_3,
 \;\; j=1,2,
\end{equation}
and from (\ref{rotmatrixident}) we obtain
\begin{eqnarray}
Q R_{\psi_2}^T R_{\varphi}^T {\bf e}_3 \hspace{-.2cm} & = & \hspace{-.2cm}
 R_{\psi_1}^T R_{\varphi}^T {\bf e}_3 \nonumber \\
R_{\varphi} R_{\psi_1}  Q R_{\psi_2}^T R_{\varphi}^T {\bf e}_3 
 \hspace{-.2cm} & = & \hspace{-.2cm}
 {\bf e}_3  \nonumber \\
R_{\varphi} R_{\psi_1} Q R_{\psi_2}^T R_{\varphi}^T
 \hspace{-.2cm} & = & \hspace{-.2cm}
 R_{\vartheta} \nonumber \\
Q \hspace{-.2cm} & = & \hspace{-.2cm}
R_{\psi_1}^T R_{\varphi}^T R_{\vartheta} R_{\varphi} R_{\psi_2}.
\end{eqnarray}
Note that $ R_{\vartheta} $, for arbitrary $ \vartheta $, is a 
rotation matrix that fixes the basis vector $ {\bf e}_3 $.

Defining $ \tilde{\bf b} $ by 
$ {\bf b} = R_{\psi_1}^T R_{\varphi}^T \tilde{\bf b} $, after
some calculation, one can show that 
\begin{equation}
\label{gedecomp}
\left[ \hspace{-.15cm} 
 \begin{array} {c c} Q & {\bf b} \\
 {\bf 0} & 1 \end{array} \hspace{-.15cm} \right] 
 \hspace{-.1cm} = \hspace{-.1cm} \left[ \hspace{-.15cm}
 \begin{array} {c c} R_{\psi_1}^T  & {\bf 0} \\
 {\bf 0} & 1 \end{array} \hspace{-.15cm} \right] \hspace{-.15cm}
\left[ \hspace{-.15cm}
 \begin{array} {c c}  R_{\varphi}^T & {\bf 0} \\
 {\bf 0} & 1 \end{array} \hspace{-.15cm} \right] \hspace{-.15cm}
\left[ \hspace{-.15cm}
 \begin{array} {c c} R_{\vartheta} & \tilde{\bf b} \\
 {\bf 0} & 1 \end{array} \hspace{-.15cm} \right] \hspace{-.15cm}
\left[ \hspace{-.15cm} \begin{array} {c c} R_{\varphi} & {\bf 0} \\
 {\bf 0} & 1 \end{array} \hspace{-.15cm} \right] \hspace{-.15cm}
\left[ \hspace{-.15cm} \begin{array} {c c}  R_{\psi_2} & {\bf 0} \\
 {\bf 0} & 1 \end{array} \hspace{-.15cm} \right] \hspace{-.1cm},
\end{equation}
\begin{equation}
\label{tildebfinal}
 \tilde{\bf b} = \left[ \hspace{-.1cm}
 \begin{array} {c} \frac{a}{a^2+w^2}\sin\vartheta \\ 
 \frac{a}{a^2+w^2}(1-\cos\vartheta) \\ \tilde{b}_3 \end{array}
 \hspace{-.1cm} \right].
\end{equation}
Thus, $ g_e $ can be decomposed as a product of five rigid motions
(four of which represent pure rotations), and contains two free
parameters - $ \vartheta $ and $ \tilde{b}_3 $ - once the control vectors
$ \Omega_1 $ and $ \Omega_2 $ are specified. 

\vspace{.25cm}

\noindent
{\bf Remark}: For purposes of interpretation of (\ref{gedecomp})
in the context of particle trajectories, we may take
$ R_{\psi_1} = R_{\psi_2} = I $, so that (\ref{gedecomp}) reduces
to
\begin{equation}
\label{gedecompinterp}
\left[ \hspace{-.1cm}
 \begin{array} {c c} Q & {\bf b} \\
 {\bf 0} & 1 \end{array} \hspace{-.1cm} \right]
 = \left[ \hspace{-.1cm}
 \begin{array} {c c}  R_{\varphi}^T & {\bf 0} \\
 {\bf 0} & 1 \end{array} \hspace{-.1cm} \right]
\left[ \hspace{-.1cm}
 \begin{array} {c c} R_{\vartheta} & \tilde{\bf b} \\
 {\bf 0} & 1 \end{array} \hspace{-.1cm} \right]
\left[ \hspace{-.1cm} \begin{array} {c c} R_{\varphi} & {\bf 0} \\
 {\bf 0} & 1 \end{array} \hspace{-.1cm} \right].
\end{equation}
To see this, recall that by definition $ g = g_1^{-1} g_2 $.
Let $ \tilde{g}_e $ be defined by
\begin{equation}
\label{getildedefn}
g_e = \left[ \hspace{-.1cm}
 \begin{array} {c c} R_{\psi_1}^T & {\bf 0} \\
 {\bf 0} & 1 \end{array} \hspace{-.1cm} \right]
 \tilde{g}_e
\left[ \hspace{-.1cm}
 \begin{array} {c c} R_{\psi_2} & {\bf 0} \\
 {\bf 0} & 1 \end{array} \hspace{-.1cm} \right].
\end{equation}
Then
\begin{eqnarray}
\label{rightmulrot}
\tilde{g}_e \hspace{-.2cm} & = & \hspace{-.2cm} \left[ \hspace{-.1cm}
 \begin{array} {c c} R_{\psi_1} & {\bf 0} \\
 {\bf 0} & 1 \end{array} \hspace{-.1cm} \right]
 g_1^{-1} g_2
\left[ \hspace{-.1cm}
 \begin{array} {c c} R_{\psi_2}^T & {\bf 0} \\
 {\bf 0} & 1 \end{array} \hspace{-.1cm} \right] \nonumber \\
\hspace{-.2cm} & = & \hspace{-.2cm} \left(g_1 \left[ \hspace{-.1cm}
 \begin{array} {c c} R_{\psi_1}^T & {\bf 0} \\
 {\bf 0} & 1 \end{array} \hspace{-.1cm} \right] \right)^{-1}
\left(g_2 \left[ \hspace{-.1cm}
 \begin{array} {c c} R_{\psi_2}^T & {\bf 0} \\
 {\bf 0} & 1 \end{array} \hspace{-.1cm} \right] \right).
\end{eqnarray} 
Thus, if we exhibit a shape equilibrium
$ \tilde{g}_e $ of the form (\ref{gedecompinterp}), we can
always write down a family of shape equilibria (\ref{getildedefn})
parameterized by $ \psi_1 $ and $ \psi_2 $, which differ only
in the orientation of the unit normal vectors of the two frames
(and are therefore indistinguishable if only the particle 
trajectories in $ \mathbb{R}^3 $ are observed).
$ \Box $ 

\vspace{.25cm}

\noindent
{\bf Proposition 1}:
Consider the two-particle system on $ G \times G $ given by
\begin{equation}
\label{se3xse3system}
\dot{g}_1 = g_1
\left[ \begin{array} {c c} \hat{\Omega}_1 & {\bf e}_1 \\ {\bf 0} & 0
 \end{array} \right], \;\;
\dot{g}_2 = g_2
\left[ \begin{array} {c c} \hat{\Omega}_2 & {\bf e}_1 \\ {\bf 0} & 0
 \end{array} \right],
\end{equation}
where $ \Omega_1 = \Omega_1(g) $, $ \Omega_2 = \Omega_2(g) $, and
$ g = g_1^{-1} g_2 $ (i.e., the controls $ \Omega_1 $ and 
$ \Omega_2 $ are arbitrary, but are $ G $-invariant).
Then there is a corresponding reduced system on $ G $
(the ``shape space'') given by 
\begin{equation}
\dot{g} = - \left[ \begin{array} {c c} \hat{\Omega}_1 & {\bf e}_1 \\ {\bf 0} & 0
 \end{array} \right] g + g
\left[ \begin{array} {c c} \hat{\Omega}_2 & {\bf e}_1 \\ {\bf 0} & 0
 \end{array} \right],
\end{equation}
(c.f. (\ref{gdot})) whose equilibria are given by (\ref{equilblockmatrix}).
Solutions of (\ref{equilblockmatrix}), with (\ref{equilblockmatrixdefn}),
require that (\ref{ctrlvectors}) hold.  
\begin{itemize}
\item[(1)] If $ w = a = 0 $, then
$ Q $ satisfies $ Q {\bf e}_1 = {\bf e}_1 $, and $ {\bf b} $ is arbitrary.
Then $ Q $ yields one free parameter, and
$ {\bf b} $ yields three free parameters. 
\item[(2)] If $ w^2 + a^2 \ne 0 $, then $ (Q,{\bf b}) $ satisfies
(\ref{gedecomp}), with $ R_{\psi_1} $, $ R_{\psi_2} $, 
$ R_{\varphi} $, and $ R_{\vartheta} $ given by (\ref{rotmatrixdefn}) 
and with $ \tilde{\bf b} $ given by (\ref{tildebfinal}).
The angle $ \varphi $ is related to $ w $ and $ a $
through (\ref{varphidefn}), and $ \vartheta $ and $ \tilde{b}_3 $ 
are free parameters.
\end{itemize}
The resulting $ (Q,{\bf b}) $ then describe the shape equilibria 
(i.e., the relative equilibria) for (\ref{se3xse3system}).

\vspace{.25cm}

\noindent
{\bf Proof}: Follows from the calculations outlined above. $ \Box $

\vspace{.25cm}

\noindent
{\bf Proposition 2}: 
Consider (\ref{se3xse3system}) as the underlying dynamics for
the evolution of two particle trajectories in $ \mathbb{R}^3 $
and their corresponding natural Frenet frames.
Then relative equilibria $ (Q,{\bf b}) $ for (\ref{se3xse3system})
correspond to the following steady-state formations of
the two particles in $ \mathbb{R}^3 $:
\begin{itemize}
\item[(1)] If $ w = a = 0 $, then the two particles move in
the same direction with arbitrary relative positions.
\item[(2)] If $ w = 0 $ but $ a \ne 0 $, then the particles
move on circular orbits with a common radius, in planes 
perpendicular to a common axis.
\item[(3)] If $ w \ne 0 $ but $ a = 0 $, then the particles
move in the same direction on collinear trajectories.
\item[(4)] If $ w \ne 0 $ and $ a \ne 0 $, then the particles
follow circular helices with the same radius, pitch, axis,
and axial direction of motion.
\end{itemize}

\vspace{.25cm}

\noindent
{\bf Proof}: Omitted due to space constraints, but follows from 
{\bf Proposition 1}, along with the {\bf Remark} 
and calculations outlined above. $ \Box $

\subsection{Shape equilibria for an n-particle system on SE(3)}

Our definition of the shape variable $ g $ for the two-particle problem
extends naturally to the
$ n $-particle problem (under the assumption that the $ n $-particle 
interaction law has $ G $ as a symmetry group).  We define
\begin{equation}
\tilde{g}_j = g_1^{-1} g_j, \;\; j=2,...,n,
\end{equation}
where $ g_1,g_2,...,g_n $ are the
group variables (each representing one of the particles), and
$ \tilde{g}_2,\tilde{g}_3,...,\tilde{g}_n $ are shape variables.
(This is analogous to the approach taken in the planar problem, 
where the corresponding group is SE(2) \cite{scltechrep,scl02,cdc03}.)  

\vspace{.25cm}

\begin{figure}
\epsfxsize=8.5cm
\epsfbox{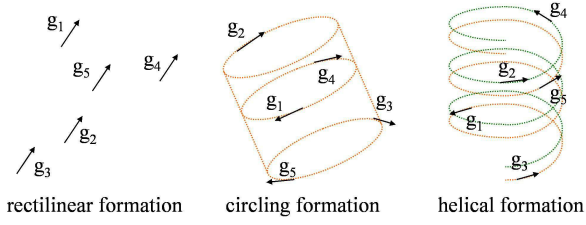}
\caption{\label{rel_eq_3d_fig} Rectilinear, circling, and helical
formations, illustrated for five particles.  The arrows represent
the unit tangent vectors to the particle trajectories.} 
\end{figure}

\noindent
{\bf Proposition 3}:
Consider 
\begin{equation}
\label{se3xse3nsystem}
\dot{g}_1 =
 g_1 \left[ \begin{array} {c c} \hat{\Omega}_1 & {\bf e}_1 \\ {\bf 0} & 0
 \end{array} \right], \;\; ..., \;\;
\dot{g}_n =
 g_n \left[ \begin{array} {c c} \hat{\Omega}_n & {\bf e}_1 \\ {\bf 0} & 0
 \end{array} \right],
\end{equation}
where $ \Omega_1,...,\Omega_n $ are $ G $-invariant controls,
as the underlying dynamics for
the evolution of $ n $ particle trajectories in $ \mathbb{R}^3. $
Then relative equilibria $ (Q_2,{\bf b}_2),...,(Q_n,{\bf b}_n) $ 
for (\ref{se3xse3nsystem})
correspond to the following steady-state formations of
the $ n $ particles in $ \mathbb{R}^3 $ (see figure \ref{rel_eq_3d_fig}):
\begin{itemize}
\item[(1)] If $ w = a = 0 $, then the $ n $ particles all move in
the same direction with arbitrary relative positions.
\item[(2)] If $ w = 0 $ but $ a \ne 0 $, then the particles
move on circular orbits with a common radius, in planes
perpendicular to a common axis.
\item[(3)] If $ w \ne 0 $ but $ a = 0 $, then the particles
move in the same direction on collinear trajectories.
\item[(4)] If $ w \ne 0 $ and $ a \ne 0 $, then the particles
follow circular helices with the same radius, pitch, axis,
and axial direction of motion.
\end{itemize}

\vspace{.25cm}

\noindent
{\bf Proof}: Omitted due to space constraints, but analogous to the
proof of {\bf Proposition 2}. $ \Box $

\vspace{.25cm}

\noindent
{\bf Remark}:
When $ w \ne 0 $ at a relative equilibrium for our model 
(\ref{se3xse3system}) of particles
evolving in $ G \times G $, the corresponding natural curvatures 
in (\ref{twouavsystem3d}) are then
in fact periodic functions of time (or arc-length parameter).  $ \Box $

\section{Rectilinear formation law}

\begin{figure}
\epsfxsize=5.5cm
\epsfbox{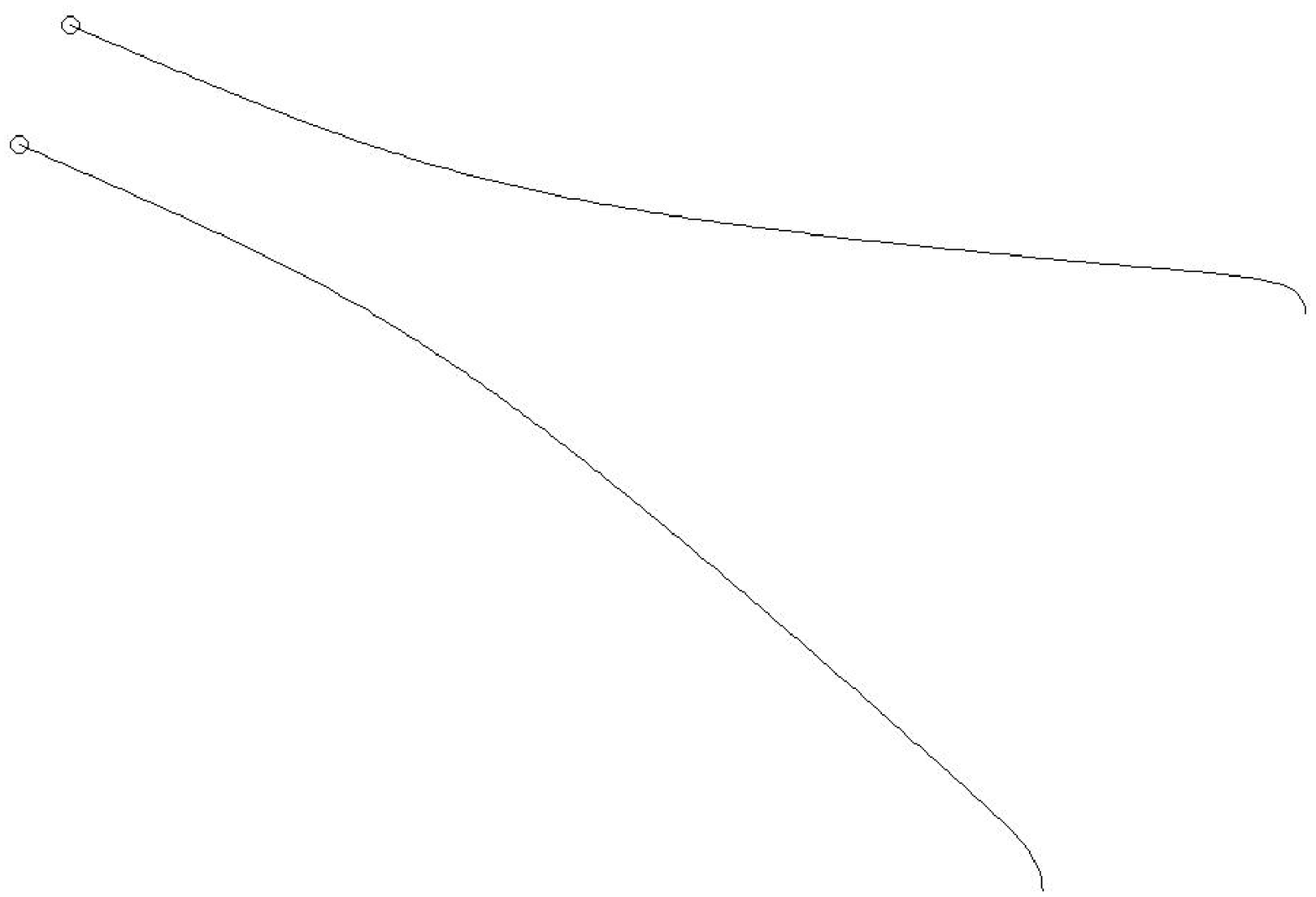}
\epsfxsize=2.8cm
\epsfbox{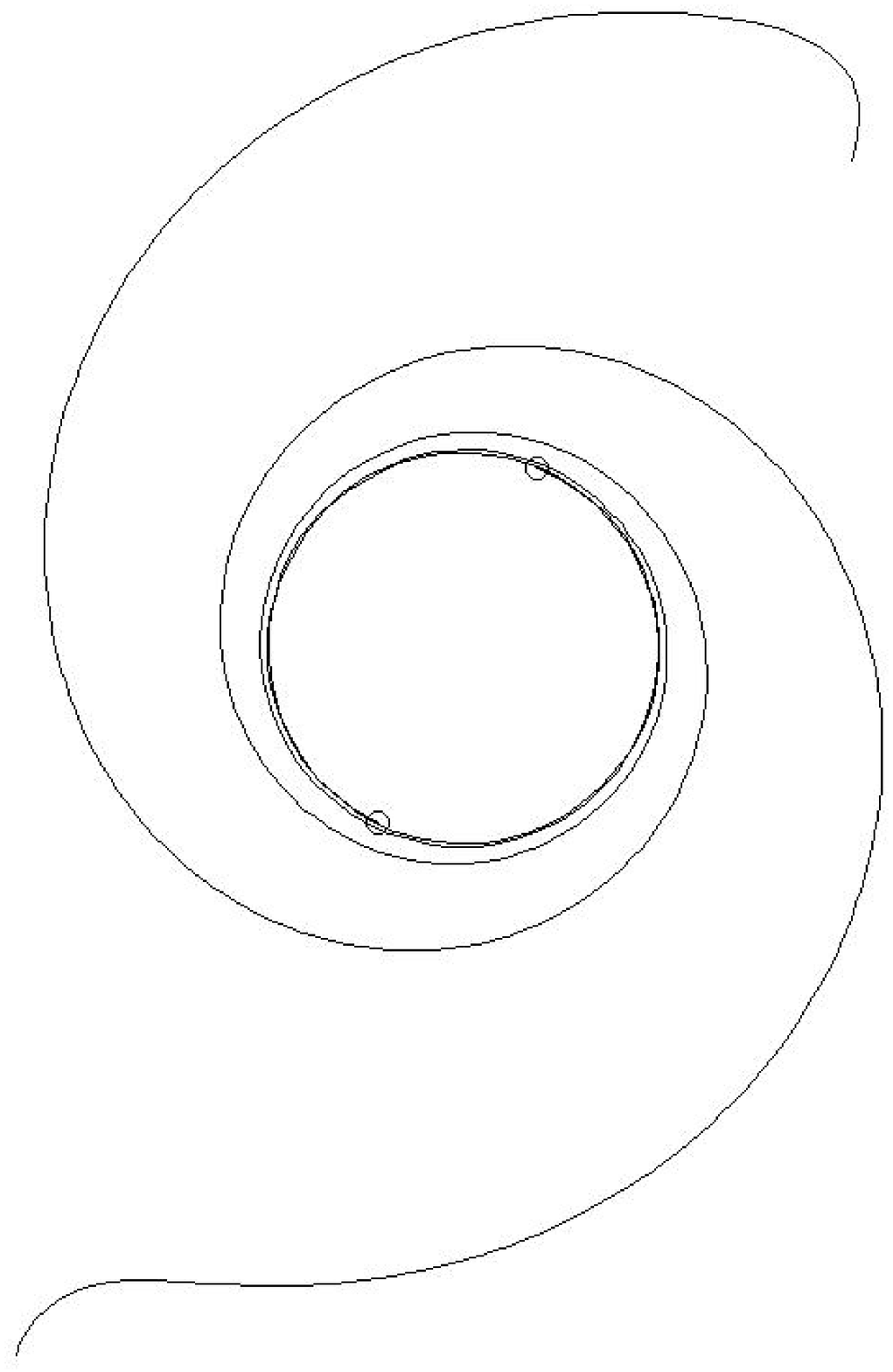}
\caption{\label{rectcirc} Convergence to a rectilinear formation (left),
and to a circling formation (right).  The trajectories, which are
three-dimensional, are viewed perpendicular to
the plane of the equilibrium formation. }
\end{figure}

The two types of equilibrium formations for which we consider
specific stabilizing control laws (for a pair of vehicles)
are rectilinear formations (in which both vehicles head in
the same direction) and circling
formations (in which both vehicles follow the same circular orbit).
Figure \ref{rectcirc} shows simulations which converge to these
two types of equilibrium formations.  
For concreteness, we use the variables $ ({\bf r}_1,{\bf x}_1,{\bf y}_1) $
and $ ({\bf r}_2,{\bf x}_2,{\bf y}_2) $,
rather than the group variables $ g_1 $ and $ g_2 $.

\begin{figure}
\hspace{1.5cm}
\epsfxsize=5cm
\epsfbox{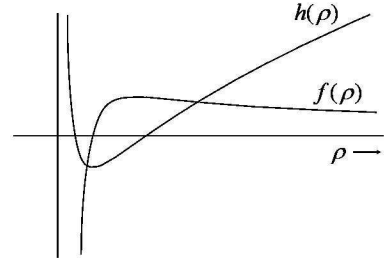}
\caption{\label{fhfig} An example of suitable functions $ f(\cdot) $
and $ h(\cdot) $ satisfying conditions (A1) and (A2) \cite{scl02}.}
\end{figure}

Consider the Lyapunov function candidate
\begin{equation}
\label{vrect}
V_{\mathit rect} = -\ln(1+{\bf x}_2\cdot {\bf x}_1) + h(|{\bf r}|),
\end{equation}
where we assume that 
\begin{itemize}
\item[$ ( \hspace{-.05cm} \mbox{A}1 \hspace{-.05cm} ) $] 
$ dh/d\rho = f(\rho) $, where $ f(\rho) $ is a Lipschitz 
continuous function on $ (0,\infty) $, so that
$ h(\rho) $ is continuously differentiable on $ (0,\infty) $;
\item[$ ( \hspace{-.05cm} \mbox{A}2 \hspace{-.05cm} ) $] 
$ \lim_{\rho \rightarrow 0} h(\rho) = \infty $,
$ \lim_{\rho \rightarrow \infty} h(\rho) = \infty $, and
$ \exists \tilde{\rho} \mbox{ such that } h(\tilde{\rho}) = 0 $.
\end{itemize}
Figure \ref{fhfig} shows an example of functions $ f(\cdot) $ and
$ h(\cdot) $ satisfying conditions (A1) and (A2).
An example of a suitable function $ f(\cdot) $ is
\begin{equation}
\label{fofr}
f(|{\bf r}|)=\alpha \left[1-\left({r_o}/{|{\bf r}|}\right)^2\right],
\end{equation}
where $ \alpha $ and $ r_o $ are positive constants.
Observe that the term
$ -\ln(1+{\bf x}_2 \cdot {\bf x}_1) $ in (\ref{vrect})
penalizes heading-direction misalignment
between the two vehicles, and the term $ h(|{\bf r}|) $
penalizes vehicle separations which are too large or too small.

Differentiating $ V_{\mathit rect} $
with respect to time along trajectories of (\ref{twouavsystem3d}) gives
\begin{eqnarray}
\label{dotvrect}
\dot{V}_{\mathit rect} \hspace{-.2cm} & = & \hspace{-.2cm}
 -\frac{\dot{\bf x}_2 \cdot {\bf x}_1 + {\bf x}_2 \cdot
 \dot{\bf x}_1}{1+{\bf x}_2\cdot{\bf x_1}}+f(|{\bf r}|)
 \frac{d}{dt}|{\bf r}| \nonumber \\
 \hspace{-.2cm} & = & \hspace{-.2cm}
 - \frac{({\bf y}_2 u_2 + {\bf z}_2 v_2)\cdot{\bf x}_1
 + {\bf x}_2 \cdot({\bf y}_1 u_1 + {\bf z}_1 v_1)}
 {1+{\bf x}_2\cdot{\bf x_1}} \nonumber \\ & & \hspace{.5cm}
 +f(|{\bf r}|) \left[
 \frac{\bf r}{|{\bf r}|}\cdot ({\bf x}_2 - {\bf x}_1)\right]
 \nonumber \\
 \hspace{-.2cm} & = & \hspace{-.2cm}
 -\frac{1}{1+{\bf x}_2\cdot{\bf x_1}} \bigg\{
 ({\bf x}_1 \cdot {\bf y}_2)u_2 +
 ({\bf x}_2 \cdot {\bf y}_1) u_1 \nonumber \\ & & \hspace{2cm} +
 ({\bf x}_1 \cdot {\bf z}_2) v_2 + ({\bf x}_2\cdot{\bf z}_1) v_1
  \nonumber \\ & & \hspace{-.2cm}
 -f(|{\bf r}|)\left(1+{\bf x}_2\cdot{\bf x_1}\right)\left[
 \frac{\bf r}{|{\bf r}|}\cdot ({\bf x}_2 - {\bf x}_1) \right] 
 \hspace{-.1cm} \bigg\}.
\end{eqnarray}
If we consider control laws of the form (\ref{twovehiclelaw3d}),
then (\ref{dotvrect}) becomes (after some calculation)
\begin{eqnarray}
\label{dotvrect2}
\dot{V}_{\mathit rect} \hspace{-.2cm} & = & \hspace{-.2cm}
 -\frac{1}{1+{\bf x}_2\cdot{\bf x_1}} \nonumber \\ 
 & & \hspace{-1.3cm} \times \hspace{-.05cm} \bigg[ \hspace{-.05cm}
 ({\bf x}_1 \cdot {\bf y}_2) F({\bf r},{\bf x}_2,{\bf y}_2,{\bf x}_1) 
 + ({\bf x}_2 \cdot {\bf y}_1) F(-{\bf r},{\bf x}_1,{\bf y}_1,{\bf x}_2)
  \nonumber \\ & & \hspace{-1.3cm} +
 ({\bf x}_1 \cdot {\bf z}_2) F({\bf r},{\bf x}_2,{\bf z}_2,{\bf x}_1)
 + ({\bf x}_2\cdot{\bf z}_1) F(-{\bf r},{\bf x}_2,{\bf z}_2,{\bf x}_1) 
 \hspace{-.05cm} \bigg].  \nonumber \\
\end{eqnarray}
It is clear from (\ref{dotvrect2}) that one choice of $ F $ which
makes $ \dot{V}_{\mathit rect} \le 0 $ is 
$ F({\bf r},{\bf x}_2,{\bf y}_2,{\bf x}_1) = \mu {\bf x}_1 \cdot {\bf y}_2, $
where $ \mu = \mu(|{\bf r}|) > 0 $.  But more generally, we consider
\begin{equation}
\label{formoff}
F({\bf r},{\bf x}_2,{\bf y}_2,{\bf x}_1) = 
 \mp \eta \left(\frac{\bf r}{|{\bf r}|}\cdot {\bf x}_2 \right)
 \left(\frac{\bf r}{|{\bf r}|}\cdot {\bf y}_2 \right) +
 \mu {\bf x}_1 \cdot {\bf y}_2,
\end{equation}
where $ \mu $ and $ \eta $ satisfy
\begin{itemize}
\item[$ ( \hspace{-.05cm} \mbox{A}3 \hspace{-.05cm} ) $] 
$ \mu(\rho) $ and $ \eta(\rho) $ are Lipschitz continuous 
on $ (0, \infty) $; 
\item[$ ( \hspace{-.05cm} \mbox{A}4 \hspace{-.05cm} ) $] 
$ \mu(|{\bf r}|) > \frac{1}{2}\eta(|{\bf r}|) > 0 $,
$ \forall |{\bf r}| \ge 0. $
\end{itemize}
(For simplicity, $ \mu $ and $ \eta $ can be taken to be constants,
rather than functions of $ |{\bf r}| $.) 

The control law given by (\ref{twovehiclelaw3d}) with (\ref{formoff})
is the natural generalization to three dimensions of the planar 
two-vehicle rectilinear law analyzed in \cite{scltechrep,scl02,cdc03}. 
As in the planar setting, we can interpret the terms involving $ f $
as steering the vehicles apart to avoid collisions (or steering them
together into formation if they are too far apart).
The terms involving $ \mu $ serve to align the vehicle headings,
and the terms involving $ \eta $ serve to align the vehicle headings
perpendicular to (or parallel to) the baseline between the vehicles.

The key to proving $ \dot{V}_{\mathit rect} \le 0 $ rests with the inequality
\begin{eqnarray}
\label{baseineq}
({\bf x}_1\cdot {\bf y}_2)\left[\frac{1}{2}({\bf x}_1\cdot {\bf y}_2)
 \mp  \left(\frac{\bf r}{|{\bf r}|}\cdot {\bf x}_2\right)
  \left(\frac{\bf r}{|{\bf r}|}\cdot {\bf y}_2\right) \right]
 \hspace{-7cm} & & \nonumber \\
 & & + ({\bf x}_2\cdot {\bf y}_1)\left[\frac{1}{2}({\bf x}_2\cdot {\bf y}_1)
 \mp   \left(\frac{\bf r}{|{\bf r}|}\cdot {\bf x}_1\right)
  \left(\frac{\bf r}{|{\bf r}|}\cdot {\bf y}_1\right) \right]
  \nonumber \\
 & & + ({\bf x}_1\cdot {\bf z}_2)\left[\frac{1}{2}({\bf x}_1\cdot {\bf z}_2)
 \mp \left(\frac{\bf r}{|{\bf r}|}\cdot {\bf x}_2\right)
  \left(\frac{\bf r}{|{\bf r}|}\cdot {\bf z}_2\right) \right]
  \nonumber \\
 & & + ({\bf x}_2\cdot {\bf z}_1)\left[\frac{1}{2}({\bf x}_2\cdot {\bf z}_1)
 \mp \left(\frac{\bf r}{|{\bf r}|}\cdot {\bf x}_1\right)
  \left(\frac{\bf r}{|{\bf r}|}\cdot {\bf z}_1\right) \right]
 \ge 0, \nonumber \\
\end{eqnarray}
which after some algebra can be shown to be equivalent to
\begin{eqnarray}
\label{baseineq2}
\left[1-({\bf x}_1 \cdot {\bf x}_2)^2\right] 
\hspace{-.3cm} & \pm & \hspace{-.3cm} 
  \bigg\{ \hspace{-.05cm} ({\bf x}_1 \cdot {\bf x}_2) \bigg[ \hspace{-.1cm}
 \left(\frac{\bf r}{|{\bf r}|}\cdot {\bf x}_1 \right)^2
 + \left(\frac{\bf r}{|{\bf r}|}\cdot {\bf x}_2 \right)^2 
 \bigg] \nonumber \\ & &  \hspace{-.2cm}
 - 2 \left(\frac{\bf r}{|{\bf r}|}\cdot {\bf x}_1 \right)
 \left(\frac{\bf r}{|{\bf r}|}\cdot {\bf x}_2 \right) \bigg\} \ge 0. 
\end{eqnarray}
If $ {\bf x}_1 = \pm {\bf x}_2 $, then (\ref{baseineq2}) holds with
equality for any choice of $ {\bf r} $.  So suppose 
$ {\bf x}_1 \ne \pm {\bf x}_2 $, and consider minimizing the expression
in (\ref{baseineq2}) over all unit vectors $ {\bf r}/|{\bf r}| $.
It is not difficult to see that (\ref{baseineq2}) achieves its minimum
for some $ {\bf r}/|{\bf r}| $ lying in the unique plane $ P $ 
containing $ {\bf x}_1 $ and $ {\bf x}_2 $ (indeed, any component of
$ {\bf r}/|{\bf r}| $ which is perpendicular to $ P $ will not contribute
to expression (\ref{baseineq2}).)  Thus, (\ref{baseineq2}) may be viewed
as a planar inequality, and we can define angle variables $ \phi_1 $ 
and $ \phi_2 $ such that 
\begin{eqnarray}
\label{anglevars}
 & & \hspace{-1cm} {\bf x}_1 \cdot {\bf x}_2  
 = \cos(\phi_2 - \phi_1), \nonumber \\
 & & \hspace{-1cm}  \left(\frac{\bf r}{|{\bf r}|}\cdot {\bf x}_1 \right)
 = \sin\phi_1, \;\;
 \left(\frac{\bf r}{|{\bf r}|}\cdot {\bf x}_2 \right)
 = \sin\phi_2.
\end{eqnarray}
After substituting (\ref{anglevars}) and applying some trigonometric 
identities, inequality (\ref{baseineq2}) becomes
\begin{equation}
\label{baseineq3}
\sin(\phi_2 - \phi_1)\left[\sin(\phi_2 - \phi_1) \pm \frac{1}{2}
 (\sin2\phi_2-\sin2\phi_1)\right] \ge 0.
\end{equation}
It can be shown that inequality (\ref{baseineq3}) does indeed hold
\cite{scltechrep}.

In the previous section, we defined shape variables
in terms of group variables in $ SE(3) $.
However, for the two-vehicle problem at hand, we can 
use the variables $ ({\bf r},{\bf x}_1,{\bf x}_2) $ 
instead, because equilibria of the 
$ ({\bf r},{\bf x}_1,{\bf x}_2) $ dynamics will include all
possible rectilinear formations.
Note that $ V_{\mathit rect} $
depends only on $ ({\bf r},{\bf x}_1,{\bf x}_2) $, as does
$ \dot{V}_{\mathit rect} $ (due to the restrictions on the control
laws we consider).  Furthermore, the 
$ ({\bf r},{\bf x}_1,{\bf x}_2) $ dynamics are self-contained
as a result of (\ref{vehctrlrestrict1})-(\ref{vehctrlrestrict2}).

\vspace{.25cm}

\noindent
{\bf Proposition 4}:
Consider the system $ ({\bf r},{\bf x}_1,{\bf x}_2) $ evolving on
$ \mathbb{R}^3 \times S^2 \times S^2 $, where $ S^2 $ is the two-sphere, 
according to (\ref{twouavsystem3d}), (\ref{twovehiclelaw3d}),
and (\ref{formoff}).  In addition, assume (A1), (A2), (A3), and (A4).
Define the set 
\begin{equation}
\Lambda = \bigg\{ ({\bf r},{\bf x}_1,{\bf x}_2) \bigg|
 {\bf x}_2 \cdot {\bf x}_1 \ne -1 \mbox{ and } |{\bf r}|>0 \bigg\}.
\end{equation}
Then any trajectory starting in $ \Lambda $ converges to the set of
equilibrium points for the $ ({\bf r},{\bf x}_1,{\bf x}_2) $-dynamics.

\vspace{.25cm}
 
\noindent
{\bf Proof}:
Observe that $ V_{\mathit rect} $ given by (\ref{vrect}) is continuously
differentiable on $ \Lambda $.  By assumption (A2) and the form of
$ V_{\mathit rect} $, we conclude that $ V_{\mathit rect} $ is radially
unbounded (i.e., $ V_{\mathit rect} \rightarrow \infty $ as
$ {\bf x}_1 \cdot {\bf x}_2 \rightarrow -1 $, as $ |{\bf r}| \rightarrow 0 $,
or as $ |{\bf r}| \rightarrow \infty $).  Therefore, for each trajectory
starting in $ \Lambda $ there exists a compact sublevel set $ \Omega $
of $ V_{\mathit rect} $ such that the trajectory remains in $ \Omega $
for all future time.  Then by LaSalle's Invariance Principle \cite{khalil},
the trajectory converges to the largest invariant set $ M $ of the
set $ E $ of all points in $ \Omega $ where $ \dot{V}_{\mathit rect}  = 0 $.
The set $ E $ in this case is the set of all points 
$ ({\bf r},{\bf x}_1,{\bf x}_2) \in \Omega $ such that 
$ {\bf x}_2 = {\bf x}_1 $.  Certainly if $ {\bf x}_1 = {\bf x}_2 
= \pm {\bf r}/|{\bf r}| $, then $ u_1 = u_2 = v_1 = v_2 = 0 $ and
the trajectory remains in $ E $ for all future time. 
Similarly, if $ {\bf r} \cdot {\bf x}_1 = {\bf r} \cdot {\bf x}_2 = 0 $
and $ f(|{\bf r}|) = 0 $, then $ u_1 = u_2 = v_1 = v_2 = 0 $ and
the trajectory remains in $ E $ for all future time.
Otherwise, we have the following expressions for the time-evolution of
the quantities $ {\bf r}\cdot{\bf x}_1 $ and $ {\bf r}\cdot {\bf x}_2 $
at points in $ E $:
\begin{eqnarray}
\label{ddtrdotx1}
\frac{d}{dt} ({\bf r}\cdot {\bf x}_1) \hspace{-1.6cm} & & \hspace{.7cm} =
 \dot{\bf r}\cdot {\bf x}_1 + {\bf r}\cdot \dot{\bf x}_1 \nonumber \\
\hspace{-.2cm} & = & \hspace{-.2cm}
 ({\bf x}_2 - {\bf x}_1)\cdot {\bf x}_1 + {\bf r}\cdot({\bf y}_1 u_1
 + {\bf z}_1 v_1) \nonumber \\
\hspace{-.2cm} & = & \hspace{-.2cm}
 ({\bf r}\cdot {\bf y}_1) \hspace{-.05cm} \left[ \hspace{-.05cm} \mp \eta
 \hspace{-.05cm}
 \left(\hspace{-.05cm} \frac{\bf r}{|{\bf r}|}\cdot {\bf x}_1
 \hspace{-.05cm}\right)
 \left(\hspace{-.05cm} \frac{\bf r}{|{\bf r}|}\cdot {\bf y}_1
 \hspace{-.05cm}\right) \hspace{-.05cm} - \hspace{-.05cm}
 f(|{\bf r}|)\left( \hspace{-.05cm} -\frac{\bf r}{|{\bf r}|}\cdot {\bf y}_1
 \hspace{-.05cm}\right) \hspace{-.05cm} \right]
 \nonumber \\
\hspace{-.2cm} & + & \hspace{-.2cm}
  ({\bf r}\cdot {\bf z}_1) \hspace{-.05cm} \left[ \hspace{-.05cm} \mp \eta
 \hspace{-.05cm}
 \left( \hspace{-.05cm} \frac{\bf r}{|{\bf r}|}\cdot {\bf x}_1
 \hspace{-.05cm} \right)
 \left( \hspace{-.05cm} \frac{\bf r}{|{\bf r}|}\cdot {\bf z}_1
 \hspace{-.05cm} \right) \hspace{-.05cm} - \hspace{-.05cm}
 f(|{\bf r}|)\left(\hspace{-.05cm} -\frac{\bf r}{|{\bf r}|}\cdot {\bf z}_1
 \hspace{-.05cm}\right) \hspace{-.05cm} \right]
 \nonumber \\ 
\hspace{-.2cm} & = & \hspace{-.2cm}
|{\bf r}|\left[1 - \left(\frac{\bf r}{|{\bf r}|}\cdot {\bf x}_1\right)^2
 \right] \left[\mp \eta\left(\frac{\bf r}{|{\bf r}|}\cdot {\bf x}_1\right)
 + f(|{\bf r}|)\right],
\end{eqnarray}
and similarly,
\begin{equation}
\label{ddtrdotx2}
\frac{d}{dt} ({\bf r}\cdot {\bf x}_2) = 
|{\bf r}| \hspace{-.1cm}
 \left[1 - \left( \frac{\bf r}{|{\bf r}|}\cdot {\bf x}_1 \right)^2
 \right] \hspace{-.1cm}
  \left[\mp \eta\left( \frac{\bf r}{|{\bf r}|}\cdot {\bf x}_1 \right)
 - f(|{\bf r}|)\right] \hspace{-.05cm} . 
\end{equation}
If $ {\bf x}_1 \ne \pm {\bf r}/|{\bf r}| $ and $ f(|{\bf r}|) \ne 0 $, then
$ \frac{d}{dt} ({\bf r}\cdot {\bf x}_1) \ne
 \frac{d}{dt} ({\bf r}\cdot {\bf x}_2), $
and it follows that the trajectory leaves $ E $.   
If $ f(|{\bf r}|) = 0 $, then the trajectory remains in $ E $, but
the only invariant subset of $ E $ with $ f(|{\bf r}|) = 0 $ 
also has $ {\bf r} \cdot {\bf x}_1 = 0 $ (or 
$ {\bf x}_1 = \pm {\bf r}/|{\bf r}| $).
Therefore, the largest invariant set contained in $ E $ may be expressed as 
\begin{eqnarray}
M \hspace{-.2cm} & = & \hspace{-.2cm}
  \bigg(\bigg\{({\bf r},{\bf x}_1,{\bf x}_2) \bigg|
 {\bf x}_1 = {\bf x}_2, \; \; {\bf r}\cdot {\bf x}_1 = 0, \;\;
 f(|{\bf r}|) = 0 \bigg\}
  \nonumber \\ & & \cup
 \bigg\{({\bf r},{\bf x}_1,{\bf x}_2) \bigg|
 {\bf x}_1 = {\bf x}_2 = \pm \frac{\bf r}{|{\bf r}|} \bigg\} \bigg)\cap
 \Omega.
\end{eqnarray}
Clearly $ M $ is contained in the set of equilibria  
of the $ ({\bf r},{\bf x}_1,{\bf x}_2) $-dynamics.  To see that
there are no other equilibria in $ \Omega $, we observe
that at equilibrium, $ \dot{\bf r} = {\bf x}_2 - {\bf x}_1 = 0 $,
and hence $ {\bf x}_2 = {\bf x}_1 $.  Since at equilibrium, we
must also have
$ \frac{d}{dt} ({\bf r}\cdot {\bf x}_1) 
= \frac{d}{dt} ({\bf r}\cdot {\bf x}_2) = 0, $
we see from equations (\ref{ddtrdotx1}) and (\ref{ddtrdotx2})
that there are no equilibria in $ \Omega $ apart from those
contained in $ M $. $ \Box $

\vspace{.25cm}

\noindent
{\bf Remark}:
If $ f $ is given by (\ref{fofr}), then $ f(|{\bf r}|) = 0 $
is equivalent to $ |{\bf r}| = r_o $.  Thus, the set of equilibria consists
of formations with both vehicles heading in the same direction, and
for one type of formation, the motion of the vehicles is
perpendicular to the baseline between them with an intervehicle 
distance equal to $ r_o $.  For the other type of formation,
both vehicles follow the same straight-line trajectory, with
one leading the other by an arbitrary distance.
The stability of these equilibria depend on the choice of parameters,
and can be further analyzed using linearization. 

\vspace{.25cm}

\noindent
{\bf Remark}:
We can express $ V_{\mathit rect} $ in terms of the group variable
$ g = g_1^{-1} g_2 $ as 
\begin{equation}
V_{\mathit rect} = -\ln(1+g_{11})+h(r),
\end{equation}
and the control law as 
\begin{eqnarray}
\label{rectctrlgroup}
u_1 \hspace{-.3cm} & = & \hspace{-.3cm}
 \mp \eta(r) \left(\frac{g_{14}g_{24}}{r^2}\right)
 \hspace{-.025cm} + \hspace{-.025cm} \mu(r) g_{21}
  \hspace{-.025cm} + \hspace{-.025cm}  f(r)\left(\frac{g_{24}}{r}\right),
\nonumber \\
u_2 \hspace{-.3cm} & = & \hspace{-.3cm}
 \mp \eta(r) \left(\frac{g^{14}g^{24}}{r^2}\right)
 \hspace{-.025cm} + \hspace{-.025cm} \mu(r) g^{21}
  \hspace{-.025cm} + \hspace{-.025cm}  f(r)\left(\frac{g^{24}}{r}\right),
\nonumber \\
v_1 \hspace{-.3cm} & = & \hspace{-.3cm}
 \mp \eta(r) \left(\frac{g_{14}g_{34}}{r^2}\right)
  \hspace{-.025cm} + \hspace{-.025cm} \mu(r) g_{31}
  \hspace{-.025cm} + \hspace{-.025cm} f(r)\left(\frac{g_{34}}{r}\right),
\nonumber \\
v_2 \hspace{-.3cm} & = & \hspace{-.3cm}
 \mp \eta(r) \left(\frac{g^{14}g^{34}}{r^2}\right)
 \hspace{-.025cm} + \hspace{-.025cm} \mu(r) g^{31}
  \hspace{-.025cm} + \hspace{-.025cm} f(r)\left(\frac{g^{34}}{r}\right),  
\end{eqnarray}
where $ g = \{g_{ij} \} $, $ g^{-1} = \{ g^{ij} \} $, and
$ r = \sqrt{g_{14}^2+g_{24}^2+g_{34}^2} $. $ \Box $

\section{Circling formation law}

Consider the Lyapunov function candidate
\begin{equation}
\label{vcircdefn}
V_{\mathit circ} = -\ln\left[1 \hspace{-.05cm} - \hspace{-.05cm}
 {\bf x}_2\cdot {\bf x}_1
 \hspace{-.05cm} + \hspace{-.05cm} 2  
 \left( \hspace{-.05cm} \frac{\bf r}{|{\bf r}|}\cdot {\bf x}_2 
 \hspace{-.05cm} \right) \hspace{-.05cm}
 \left( \hspace{-.05cm} \frac{\bf r}{|{\bf r}|}\cdot {\bf x}_1 \hspace{-.05cm}
 \right) \right] + h(|{\bf r}|),
\end{equation}
where we assume 
\begin{itemize}
\item[$( \hspace{-.05cm} \mbox{A} \hspace{-.05cm} 1 \hspace{-.05cm} \mbox{'}
 \hspace{-.05cm})$] 
$ dh/d\rho = f(\rho)-2/\rho $, where $ f(\rho) $ is a Lipschitz
continuous function on $ (0,\infty) $, so that
$ h(\rho) $ is continuously differentiable on $ (0,\infty) $;
\end{itemize}
and (A2).  It can be shown that
\begin{equation}
\label{vcirclnterm}
1 - {\bf x}_2\cdot {\bf x}_1 + 2
 \left( \frac{\bf r}{|{\bf r}|}\cdot {\bf x}_2 \right) 
 \left( \frac{\bf r}{|{\bf r}|}\cdot {\bf x}_1 \right) \ge 0,
\end{equation}  
and the function $ f $ given by (\ref{fofr}) can be used here, as well.
The term $ h(|{\bf r}|) $ in (\ref{vcircdefn}) penalizes vehicle
separations which are two large or too small.  The natural-log term
in (\ref{vcircdefn}) involves the relative headings of the vehicles,
as well as the relative orientations of the headings with respect to
the baseline between the vehicles.

Differentiating $ V_{\mathit circ} $ along trajectories of 
(\ref{twouavsystem3d}) 
and plugging in (\ref{twovehiclelaw3d}) gives
\begin{eqnarray}
\dot{V}_{\mathit circ} \hspace{-.25cm} & = & \hspace{-.25cm}
 -\frac{1}{1-{\bf x}_2 \cdot {\bf x}_1 +
 2\left(\frac{\bf r}{|{\bf r}|}\cdot {\bf x}_2 \right)
 \left(\frac{\bf r}{|{\bf r}|}\cdot {\bf x}_1 \right)} \nonumber \\
 & & \hspace{-1.3cm} \times \hspace{-.05cm}
  \Bigg\{ \hspace{-.1cm} \left[ \hspace{-.05cm} -{\bf x}_1 \cdot {\bf y}_2
 \hspace{-.05cm} + \hspace{-.05cm} 2 \hspace{-.05cm} \left( \hspace{-.05cm}
 \frac{\bf r}{|{\bf r}|}\cdot {\bf x}_1 \hspace{-.05cm} \right) \hspace{-.05cm}
 \left(\hspace{-.05cm} \frac{\bf r}{|{\bf r}|}\cdot {\bf y}_2 \hspace{-.05cm}
 \right) \hspace{-.05cm} \right] \hspace{-.05cm}
 F({\bf r},{\bf x}_2,{\bf y}_2,{\bf x}_1)
 \nonumber \\ & & \hspace{-1.1cm}
 + \hspace{-.05cm} \left[ \hspace{-.05cm} -{\bf x}_2 \cdot {\bf y}_1
 \hspace{-.05cm} + \hspace{-.05cm} 2 \hspace{-.05cm}
 \left( \hspace{-.05cm} \frac{\bf r}{|{\bf r}|}\cdot {\bf x}_2 \hspace{-.05cm}
  \right) \hspace{-.05cm}
 \left(\hspace{-.05cm} \frac{\bf r}{|{\bf r}|}\cdot {\bf y}_1 \hspace{-.05cm}
 \right) \hspace{-.05cm} \right] \hspace{-.05cm}
 F(-{\bf r},{\bf x}_1,{\bf y}_1,{\bf x}_2)
 \nonumber \\ & & \hspace{-1.1cm}
 + \hspace{-.05cm} \left[ \hspace{-.05cm} -{\bf x}_1 \cdot {\bf z}_2
 \hspace{-.05cm} + \hspace{-.05cm} 2 \hspace{-.05cm}
 \left(\hspace{-.05cm} \frac{\bf r}{|{\bf r}|}\cdot {\bf x}_1 \hspace{-.05cm}
 \right) \hspace{-.05cm}
 \left( \hspace{-.05cm} \frac{\bf r}{|{\bf r}|}\cdot {\bf z}_2 \hspace{-.05cm}
  \right) \hspace{-.05cm} \right] \hspace{-.05cm}
 F({\bf r},{\bf x}_2,{\bf z}_2,{\bf x}_1)
 \nonumber \\ & & \hspace{-1.1cm}
 + \hspace{-.05cm} \left[ \hspace{-.05cm} -{\bf x}_2 \cdot {\bf z}_1
 \hspace{-.05cm} + \hspace{-.05cm} 2 \hspace{-.05cm}
 \left( \hspace{-.05cm} \frac{\bf r}{|{\bf r}|}\cdot {\bf x}_2 \hspace{-.05cm}
 \right) \hspace{-.05cm}
 \left( \hspace{-.05cm} \frac{\bf r}{|{\bf r}|}\cdot {\bf z}_1 \hspace{-.05cm}
 \right) \hspace{-.05cm} \right] \hspace{-.05cm}
 F(-{\bf r},{\bf x}_1,{\bf z}_1,{\bf x}_2)
 \hspace{-.05cm} \Bigg\} \hspace{-.02cm}. \nonumber \\
\end{eqnarray}
In place of (\ref{formoff}), we use
\begin{eqnarray}
\label{formoffcirc}
F({\bf r},{\bf x}_2,{\bf y}_2,{\bf x}_1) 
\hspace{-.2cm} & = & \hspace{-.2cm}
 \pm \eta \left(\frac{\bf r}{|{\bf r}|}\cdot {\bf x}_2\right)
  \left(\frac{\bf r}{|{\bf r}|}\cdot {\bf y}_2\right)
 \nonumber \\ & & \hspace{-2cm} + \mu \hspace{-.05cm}  \left[ \hspace{-.05cm}
 -{\bf x}_1 \cdot {\bf y}_2 
 + 2\left(\frac{\bf r}{|{\bf r}|}\cdot {\bf x}_1 \right)
 \left(\frac{\bf r}{|{\bf r}|}\cdot {\bf y}_2 \right) \hspace{-.05cm} \right] 
 \hspace{-.1cm}, 
\end{eqnarray}
where we assume (A3) and (A4).

The key to proving $ \dot{V}_{\mathit circ} \le 0 $ can then be shown to
rest with the inequality
\begin{eqnarray}
\label{circineqxr}
 1 - \left[-{\bf x}_2 \cdot {\bf x}_1
 +2\left(\frac{\bf r}{|{\bf r}|}\cdot {\bf x}_2 \right)
 \left(\frac{\bf r}{|{\bf r}|}\cdot {\bf x}_1 \right)\right]^2
\hspace{-6.2cm} & & \nonumber \\
 \hspace{-.2cm} & \pm & \hspace{-.2cm}
 \Bigg\{{\bf x}_2 \cdot {\bf x}_1 + \left[-{\bf x}_2 \cdot {\bf x}_1
 + 2\left(\frac{\bf r}{|{\bf r}|}\cdot {\bf x}_1\right)
 \left(\frac{\bf r}{|{\bf r}|}\cdot {\bf x}_2 \right) \right]
 \nonumber \\ & & \hspace{1.5cm} \times
 \left[1-\left(\frac{\bf r}{|{\bf r}|}\cdot {\bf x}_1\right)^2
 - \left(\frac{\bf r}{|{\bf r}|}\cdot {\bf x}_2 \right)^2 \right] \Bigg\}
 \nonumber \\
 \hspace{-.2cm} & \ge & \hspace{-.2cm}
 0.
\end{eqnarray}
Using a similar technique as was used above to pass from inequality
(\ref{baseineq2}) to inequality (\ref{baseineq3}), we can show that
(\ref{circineqxr}) also becomes (essentially) inequality (\ref{baseineq3}).

\vspace{.25cm}

\noindent
{\bf Proposition 5}:
Consider the system $ ({\bf r},{\bf x}_1,{\bf x}_2) $ evolving on
$ \mathbb{R}^3 \times S^2 \times S^2 $, 
according to (\ref{twouavsystem3d}), (\ref{twovehiclelaw3d}),
and (\ref{formoffcirc}).  In addition, assume (A1'), (A2), (A3), and (A4).
Define the set
\begin{eqnarray}
\Lambda' \hspace{-.25cm} & = & \hspace{-.25cm}
  \bigg\{ ({\bf r},{\bf x}_1,{\bf x}_2) \bigg|
1 \hspace{-.05cm} - \hspace{-.05cm}
 {\bf x}_2\cdot {\bf x}_1
 \hspace{-.05cm} + \hspace{-.05cm} 2
 \left( \hspace{-.05cm} \frac{\bf r}{|{\bf r}|}\cdot {\bf x}_2
 \hspace{-.05cm} \right) \hspace{-.05cm}
 \left( \hspace{-.05cm} \frac{\bf r}{|{\bf r}|}\cdot {\bf x}_1 \hspace{-.05cm}
 \right) \hspace{-.05cm} \ne \hspace{-.05cm} 0 
 \nonumber \\ & & \hspace{1.5cm}
 \mbox{ and } |{\bf r}|>0 \bigg\}.
\end{eqnarray}
Then any trajectory starting in $ \Lambda' $ converges to the set 
\begin{eqnarray}
\tilde{M}' \hspace{-.3cm} & = & \hspace{-.3cm} 
 \bigg( \hspace{-.05cm} \bigg\{ \hspace{-.05cm} 
 ({\bf r},{\bf x}_1,{\bf x}_2) \bigg|
 {\bf x}_1 \hspace{-.025cm} = \hspace{-.025cm} -{\bf x}_2, \; 
 {\bf r}\cdot {\bf x}_1 \hspace{-.025cm} = \hspace{-.025cm} 0, \;
 f(|{\bf r}|)  \hspace{-.025cm} = \hspace{-.025cm} \frac{2}{|{\bf r}|} \bigg\}
  \nonumber \\ & & \cup
 \bigg\{({\bf r},{\bf x}_1,{\bf x}_2) \bigg|
 {\bf x}_1 = {\bf x}_2 = \pm \frac{\bf r}{|{\bf r}|} \bigg\} \bigg)
 \cap \Lambda'.
\end{eqnarray}
Note that elements of $ \tilde{M}' $ with $ {\bf x}_1 = -{\bf x}_2 $
correspond to the two vehicles following the same circular orbit,
separated by the diameter of the orbit, which is prescribed by the
function $ f $.  Elements of $ \tilde{M}' $ with $ {\bf x}_1 = {\bf x}_2 $
correspond to rectilinear formations in which one vehicle
leads the other by an arbitrary distance.

\vspace{.25cm}

\noindent
{\bf Proof}:
Omitted due to space constraints, but a similar approach is used as
in the proof of {\bf Proposition 4}. $ \Box $

\vspace{.25cm}

\noindent
{\bf Remark}:
We can express $ V_{\mathit circ} $ in terms of the group variables
as
\begin{equation}
V_{\mathit circ} = -\ln\left(1-g_{11}-2\frac{g_{14}g^{14}}{r^2}\right)+h(r),
\end{equation}
and the control law for circling
can also be expressed in terms of the group variables,
analogously to (\ref{rectctrlgroup}). $ \Box $

\section{Multi-vehicle formations}

\begin{figure}
\hspace{.5cm}
\epsfxsize=7.5cm
\epsfbox{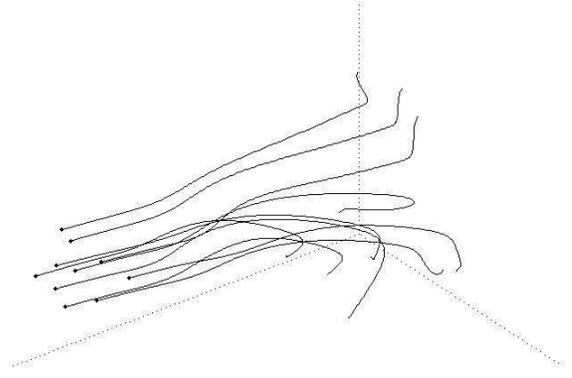}
\caption{\label{multiveh} Simulation results for ten vehicles using
generalization (\ref{multivehctrl}) of the two-vehicle rectilinear
formation control law (\ref{twovehiclelaw3d}) with (\ref{formoff})
and (\ref{fofr}).} 
\end{figure}

\begin{figure}
\epsfxsize=4cm
\epsfbox{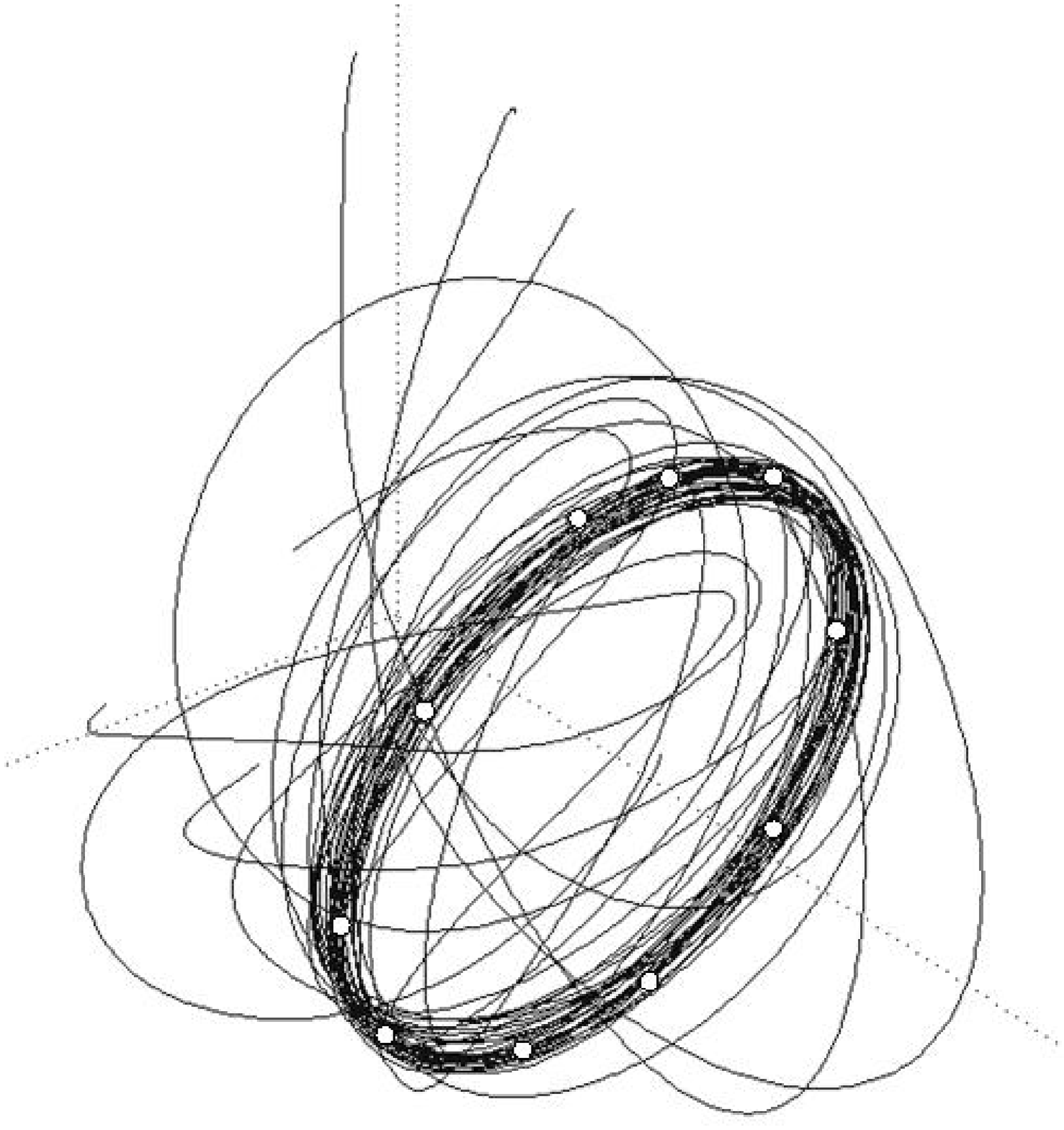}
\epsfxsize=4cm
\epsfbox{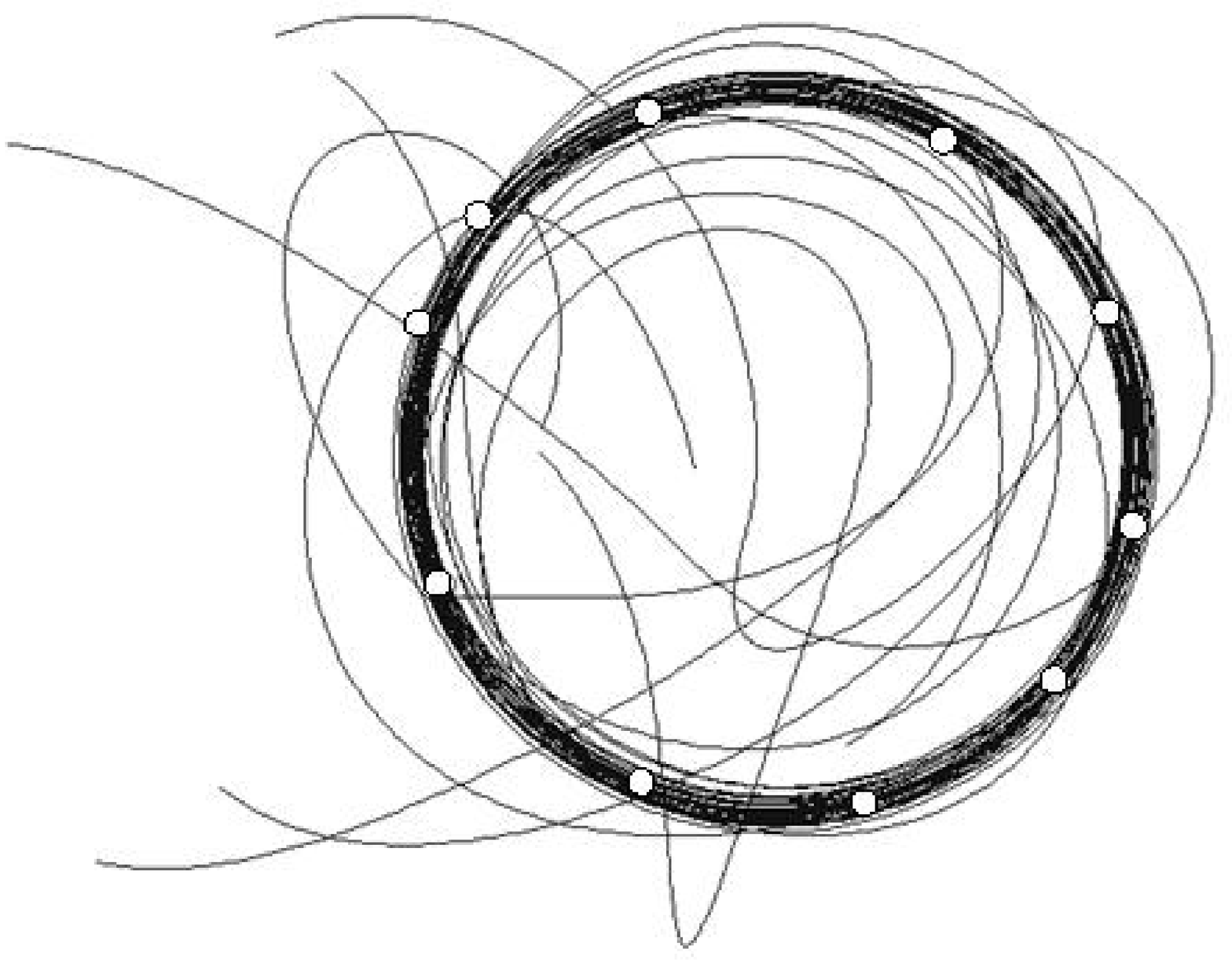}
\caption{\label{multivehcir} Simulation results for ten vehicles using
generalization (\ref{multivehctrl}) of the two-vehicle circling 
formation control law (\ref{twovehiclelaw3d}) with (\ref{formoffcirc})
and (\ref{fofr}). (The same simulation results are viewed from two 
different angles.)}
\end{figure}

One way to generalize the two-vehicle laws discussed above to $ n $ vehicles  
is to use an average of the pairwise interaction terms
used for the two-vehicle problem \cite{scltechrep,scl02,cdc03}, i.e.,
\begin{eqnarray}
\label{multivehctrl}
u_j  \hspace{-.2cm} & = & \hspace{-.2cm}
 \frac{1}{n} \sum_{k\ne j} \bigg[
 F({\bf r}_j-{\bf r}_k,{\bf x}_j,{\bf y}_j,{\bf x}_k) \nonumber \\ & &
 \hspace{1cm} - f(|{\bf r}_j-{\bf r}_k|)\left(\frac{ {\bf r}_j-{\bf r}_k }
 { |{\bf r}_j-{\bf r}_k| } \cdot {\bf y}_j \right) \bigg], \nonumber \\
v_j \hspace{-.2cm} & = & \hspace{-.2cm}
 \frac{1}{n} \sum_{k\ne j} \bigg[
 F({\bf r}_j-{\bf r}_k,{\bf x}_j,{\bf z}_j,{\bf x}_k) \nonumber \\ & & 
 \hspace{1cm} - f(|{\bf r}_j-{\bf r}_k|)\left(\frac{ {\bf r}_j-{\bf r}_k }
 { |{\bf r}_j-{\bf r}_k| } \cdot {\bf z}_j \right) \bigg], 
\end{eqnarray}
$ j=1,...,n $.  In (\ref{multivehctrl}), $ {\bf r}_j $ is the position of
the $ j^{\mbox{th}} $ vehicle, $ ({\bf x}_j,{\bf y}_j,{\bf z}_j) $ is
the corresponding natural Frenet frame, and $ (u_j,v_j) $ are
the associated natural curvatures.
Figures \ref{multiveh} and \ref{multivehcir} show simulation results
for multi-vehicle interactions of this type.  Their analysis is a 
topic of ongoing research.

\section{Acknowledgements}

This research was supported in part by the Naval Research Laboratory under
Grants No.~N00173-02-1G002, N00173-03-1G001, N00173-03-1G019, and
N00173-04-1G014; by the 
Air Force Office of Scientific Research under AFOSR Grants
No.~F49620-01-0415 and FA95500410130; by the Army Research Office
under ODDR\&E MURI01 Program Grant No.~DAAD19-01-1-0465 to the Center for
Communicating Networked Control Systems (through Boston University);
and by NIH-NIBIB grant
1 R01 EB004750-01, as part of the NSF/NIH Collaborative
Research in Computational Neuroscience Program.

\end{document}